\documentclass[11pt,leqno]{article}
\usepackage{amsmath, amsthm, amsfonts, amssymb}
\usepackage{mathrsfs}
\newfam\msbfam
\font\tenmsb=msbm10    \textfont\msbfam=\tenmsb \font\sevenmsb=msbm7
\scriptfont\msbfam=\sevenmsb \font\fivemsb=msbm5
\scriptscriptfont\msbfam=\fivemsb

\newfam\bigfam
\font\tenbig=msbm10 scaled \magstep2   \textfont\bigfam=\tenbig
\font\sevenbig=msbm7 scaled \magstep2 \scriptfont\bigfam=\sevenbig
\font\fivebig=msbm5 scaled \magstep2
\scriptscriptfont\bigfam=\fivebig

\def\dint{\displaystyle\int}

\def\dfrac{\displaystyle\frac}

\begin{document}


\begin{center} {\Large\bf  Boundedness of  commutators generated by m-th  Calder\'{o}n-Zygmund type
singular integrals and local Campanato functions on generalized  local Morrey spaces}

\vspace{0.3cm}

\noindent{\large\bf HuiXia MO\footnote[1]{
Correspondence:huixiamo@bupt.edu.cn\\
2010 AMS Mathematics Subject
Classification: 42B20, 42B25}, Hongyang XUE}\\
School of Science, Beijing University of Posts and Telecommunications, Beijing, 100876, China
\end{center}

\parskip 12pt

\noindent{\bf Abstract}\;\;{\footnotesize
Let $T_m$ be  the $m$-th Calder\'{o}n-Zygmund type singular integral.  In the paper, we consider the boundedness of
$T_m$ on the generalized product local Morrey spaces
$LM_{p_1,\varphi_1}^{\{x_0\}}\times LM_{p_2,\varphi_2}^{\{x_0\}}\times\dots\times LM_{p_m,\varphi_m}^{\{x_0\}}.$ And, the boundedness of the commutators of $T_m$  with local Campanato functions is obtained, also.}

 \vspace{0.3cm}
{\bf Key words}\;\; m-th Calder\'{o}n-Zygmud type  singular integral,
commutator, local Campanato function, generalized local Morrey space


\section{\bf Introduction}

In recent years, the multilinear singular integrals have been
attracting  attention and great developments have been achieved (see
[1-11]). The study for the multilinear singular integrals is motivated
not only by a mere quest to generalize the theory of linear
operators but also by their natural appearance in analysis.

  Meanwhile,  the commutators generated by the multilinear singular integral
and BMO functions or Lipschitz functions also attract much
attention, since the commutator is more singular than the singular
integral operator itself.

 Moreover, the classical Morrey space $M_{p,\lambda}$ were first introduced by Morrey in \cite{M} to study the local behavior of solutions to second order elliptic partial differential equations. In \cite{LL}, the authors studied the boundedess of the  multilinear Calder\'{o}n-Zygmund singular integral on the classical Morrey space $M_{p,\lambda}.$  And, in \cite{BG}, the authors introduced the local generalized Morrey space $LM_{p,\varphi}^{\{x_{0}\}}$, and they also studied the boundedness of the homogeneous singular integrals with rough kernel on these spaces.

Motivated by the works of \cite{LL,BG}, we are going to
consider the boundedness of the multilinear Calder\'{o}n-Zygmund singular integral and its commutator on the local generalized Morrey space $LM_{p,\varphi}^{\{x_{0}\}}$

Now, let us give some related notations.

We are going to be working in $\mathbb{R}^{n}$. Let $m\in\mathbb{N}$
and $K(y_{0},y_{1},\dots,y_{m})$ be a function defined away from the
diagonal $y_{0}=y_{1}=,\dots,=y_{m}$ in $(\mathbb{R}^{n})^{m+1}.$
Let $T_m$ be a multilinear operator which was initially defined on
the m-fold product of Schwartz space $\mathcal{S}(\mathbb{R}^n)$ and
take its values in the space of tempered distributions
$\mathcal{S'}(\mathbb{R}^n)$ and such that for $K$, the integral
representation below is valid:
$$\begin{array}{cl}T_m(\vec{f})(x)=T_m(f_{1},\dots,f_{m})(x) =\dint_{(\mathbb{R}^{n})^{m}}K(x,y_{1},\dots ,y_{m})
f_{1}(y_1) \dots f_{m}(y_m)dy_1\dots dy_m,\end{array}\eqno{(1.1)}$$ whenever $f_{i},$ $i=1,
\dots ,m,$ are smooth functions with compact support and
$x\notin\cap_{i=1}^{m}\mbox{supp}f_{i}.$

 Moreover, if the kernel $K$ satisfies the following size and smoothness estimates:
$$\begin{array}{cl}
|K(y_0,y_{1},\dots,y_{m})|\leq
\dfrac{C}{(\sum^{m}_{k,l=0}|y_{k}-y_{l}|)^{mn}},\end{array}\eqno{(1.2)}$$
for all $(y_0,y_{1},\dots,y_{m})\in(\mathbb{R}^{n})^{m+1}$ away from
the diagonal;
$$\begin{array}{cl}|K(y_0,\dots,y_i,\dots,y_{m})-K(y_0,\dots,y'_i,\dots,y_{m})|\leq\dfrac{C|y_i-y'_i|^{\epsilon}}
{(\sum^{m}_{k,l=0}|y_{k}-y_{l}|)^{mn+\epsilon}},\end{array}\eqno{(1.3)}$$
for some $C>0$ and $\epsilon>0,$ whenever $0\leq j\leq m$ and
$|y_i-y'_i|\leq 1/2\max_{0\leq k\leq m}|y_i-y_{k}|,$ then the kernel
is called a m-th Calder\'{o}n-Zygmund kernel and the collection of
such functions is denoted by $m-CZK(C,\epsilon)$. Let $T_m$ be as in
(1.1) with a $m-CZK(C,\epsilon)$ kernel, then $T_m$ is called a m-th
Calder\'{o}n-Zygmund type singular integral  and  the collection of
these operators is denoted by $m-CZO$.

Now, we define the  commutators generated by the m-th multilinear Calder\'{o}n-Zygmund type singular integral as follows.

Let $\vec{b}=(b_{1},\dots,b_{m})$ be a finite family of
locally integrable functions, then the commutators generated by the m-th
Calder\'{o}n-Zygmund type singular integral and $\vec{b}$ is
defined by:
$$T_{m}^{\vec{b}}(\vec{f})(x) =\dint_{(\mathbb{R}^{n})^{m}}K(x,y_{1},\dots ,y_{m})
\prod\limits_{i=1}^{m}(b_{i}(x)-b_{i}(y_{i}))f_{i}(y_i)dy_1\dots dy_m.$$

In the following,  we will establish the boundedness of $T_m$ on generalized product local Morrey spaces. And,  we also consider the boundedness of the  commutators generated by the m-th Calder\'{o}n-Zygmund type singular integral $T_m$
and the local Campanato function on generalized product local Morrey spaces.

\section{\bf Some notations and  lemmas}

\textbf{Definition 2.1}\cite{BG}~Let $\varphi(x,r)$ be a positive measurable function on $\mathbb{R}^n\times(0,\infty)$ and $1\leq p\leq\infty.$ For any fixed $\emph{x}_{0}\in\mathbb{R}^n,$ a function $f\in L_{loc}^{q}$ is said to belong to the local Morrey space, if
$$\|f\|_{LM^{\{x_0\}}_{p,\varphi}}=\sup\limits_{r>0}\varphi^{-1}(x_0,r)|B(x_0,r)|^{-\frac{1}{p}}\|f\|_{\emph{L}_p(B(x_0,r))}<\infty.$$
And, we denote
$$\emph{LM}^{\{x_0\}}_{p,\varphi}\equiv\emph{LM}^{\{x_0\}}_{p,\varphi}(\mathbb{R}^n)=\{f\in L_{loc}^{q}(\mathbb{R}^n):\|f\|_{{LM}^{\{x_0\}}_{p,\varphi}}<\infty\}.$$

According to this definition, we recover the local Morrey space$ \emph{LM}^{\{x_0\}}_{p,\lambda}$ under the choice $\varphi(x_0,r)=r^{\frac{\lambda-n}{p}}.$

\textbf{Definition 2.2}\cite{BG}~Let $1\leq q<\infty$ and $0\leq\lambda<1/n.$ A function $f\in L_{loc}^{q}(\mathbb{R}^n)$ is said to belong to the space $LC_{q,\lambda}^{\{x_0\}}$ (local Campanato space), if
$$\|f\|_{LC_{q,\lambda}^{\{x_0\}}}=\sup\limits_{r>0}\biggl(\dfrac{1}{|B(x_0,r)|^{1+\lambda q}}\dint_{B(x_,r)}|f(y)-f_{B(x_0,r)}|^{q}dy\biggr)^{1/q}<\infty,$$
where $$f_{B(x_0,r)}=\dfrac{1}{|B(x_0,r)|}\dint_{B(x_0, r)}f(y)dy.$$

Define $$LC_{q,\lambda}^{\{x_0\}}(\mathbb{R}^n)=\{f\in L_{loc}^{q}(\mathbb{R}^n):\|f\|_{LC_{q,\lambda}^{\{x_0\}}}<\infty\}.$$

{\bf Remark.}\cite{BG} Note that, the central $BMO$ space $CBMO_{q}(\mathbb{R}^n)=LC_{q,0}^{\{0\}}(\mathbb{R}^n),$
 $CBMO_{q}^{\{x_0\}}(\mathbb{R}^n)=LC_{q,0}^{\{x_0\}}(\mathbb{R}^n),$
and  $BMO_{q}(\mathbb{R}^n)\subset \bigcap_{q>1}CBMO_{q}^{\{x_0\}}(\mathbb{R}^n).$ Moreover, one can imagine that the behavior of  $CBMO_{q}^{\{x_0\}}(\mathbb{R}^n)$ may be quite different from that of $BMO(\mathbb{R}^n),$ since there is no analogy of the John-Nirenberg inequality of $BMO$ for the space $CBMO_{q}^{\{x_0\}}(\mathbb{R}^n).$

{\bf Lemma 2.1}~Let $1<q<\infty,$ $0<r_2<r_1$ and  $b\in LC_{q,\lambda}^{\{x_0\}},$ then
$$\begin{array}{cl}\bigg(\dfrac{1}{|B(x_0,r_1)|^{1+\lambda q}}\dint_{B(x_0,r_1)}|b(x)-b_{B(x_0,r_2)}|^{q}dx\biggr)^{1/q}\leq C\biggl(1+\ln\dfrac{r_1}{r_2}\biggr)\|b\|_{LC_{q,\lambda}^{\{x_0\}}}.\end{array}\eqno{(2.1)}$$
And, from this inequality, we have
$$\begin{array}{cl}|b_{B(x_0,r_1)}-b_{B(x_0,r_2)}|\leq C\biggl(1+\ln\dfrac{r_1}{r_2}\biggr)|B(x_0,r_1)|^{\lambda}\|b\|_{LC_{q,\lambda}^{\{x_0\}}}.\end{array}\eqno{(2.2)}$$

In this section, we are going to use the following statement on the boundedness of the weighted Hardy operator:
$$\emph{H}_{w}g(t):=\int^\infty_{t}g(s)w(s)ds,~0<t<\infty,$$
where $\emph{w}~$is a fixed function non-negative and measurable on $(0,\infty).$

{\bf Lemma 2.2}\cite{G1,G2}~Let $v_1,v_2$ and $w$ be positive almost everywhere and measurable functions on $(0,\infty).$ The inequality\\
$$\begin{array}{cl}ess ~\sup\limits_{t>0}v(2t)\emph{H}_wg(t)\leq Cess ~\sup\limits_{t>0}v_1(t)g(t)\end{array}\eqno{(2.3)}$$
holds for some $C>0$ and all non-negative and non-decreasing g on$~(0,\infty)~$if and only if
$$B:ess \sup\limits_{t>0}v_2(t)\int^\infty_{t}\frac{w(s)ds}{ess ~\sup_{s<\tau<\infty} v_1(\tau)}ds<\infty.$$ Moreover, if $\tilde{C}$ is the minimum value of $C$ in (2.3), then $\tilde{C}=B$.

{\bf Lemma 2.3}\cite{KS}~\;\;Let $T_m$ be a $m-CZO$.
 Suppose that $1\leq  p_{1},\cdots,p_{m}<\infty$ and $1/p=1/p_{1}+\cdots +1/p_{m}.$
 If $p_{i}>1,i=1,\cdots,m,$ then there exists a constant $C>0,$ such that
  $$\|T_{m}\vec{f}\|_{L^{p}}\leq C\prod\limits_{i=1}^{m}\|f_{i}\|_{L^{p_{i}}}.$$

\section{ M-th Calder\'{o}n-Zygmund type  singular integral operator on generalized product local Morrey space}

{\bf Theorem 3.1}\;\ Let $x_0\in{\mathbb{R}^n},$ $1<p, p_1,p_2,\dots,p_m<\infty,$ such that $1/p=1/p_1+1/p_2+\dots+p_m.$ Then
the inequality
$$\|T_m(\vec{f})\|_{L^p(B(x_0,r))}\lesssim r^{n/p}\dint^{\infty}_{2r}\prod\limits_{i=1}^{m}\|f_i\|_{L^{p_i}(B(x_0,r))}t^{-n/p-1}dt$$
holds for any ball $B(x_0,r)$ and  all $f_i\in L^{p_i}_{loc}(\mathbb{R}^n),$ $i=1,2,\dots,m.$

{\bf Proof.}  Without loss of generality, it is suffice to show that the conclusion holds for
 $T_{2}(f_{1},f_{2}).$

Let $B=B(x_0,r).$  And, we write $f_1=f^0_1+f^\infty_1$ and $f_2=f^0_2+f^\infty_2,$ where $f^0_i=f_i\chi_{2B},$ $f^\infty_i=f_i\chi_{{(2B)}^c},$ for $i=1.2.$ Thus, we have
$$\begin{array}{cl}
&\|T_2(f_1,f_2)\|_{L^p(B(x_0,r))}\\
\leq&\|T_2(f^0_1,f^0_2)\|_{L^p(B)}+\|T_2(f^0_1,f^\infty_2)\|_{L^p(B)}+\|T_2(f^\infty_1,f^0_2)\|_{L^p(B)}+\|T_2(f^\infty_1,f^\infty_2)\|_{L^p(B)}\\
=:&I+II+III+IV.
\end{array}$$

Using the $L^p$ boundedness of $T_2$(Lemma 2.3), we have
$$\begin{array}{cl}
  I  \lesssim&\|f_1\|_{L^{p_1}(2B)}\|f_2\|_{L^{p_2}(2B)} \\
 \lesssim&r^{\frac{n}{p}}\|f_1\|_{L^{p_2}(2B)}\|f_2\|_{L^{p_2}(2B)}\dint^\infty_{2r}\frac{dt}{t^{\frac{n}{p}+1}}\\
  \leq&r^\frac{n}{p}\dint^\infty_{2r}\|f_1\|_{L^{p_1}(B(x_0,t))}\|f_2\|_{L^{p_2}(B(x_0,t))}\frac{dt}{t^{\frac{n}{p}+1}}.\\
\end{array}\eqno{(3.1)}$$

Moreover, when $x\in B(x_0, r)$ and $y\in{(2B)}^c,$ we have $$\begin{array}{cl}\dfrac{1}{2}|x_0-y|\leq|x-y|\leq\dfrac{3}{2}|x_0-y|.\end{array}\eqno{(3.2)}$$

Then, it follows from (1.2) that
$$\begin{array}{cl}
|T_2(f^0_1,f^\infty_2)(x)|
  &\lesssim\dint_{\mathbb{R}^n}\dint_{\mathbb{R}^n}\dfrac{|f^0_1(y_1)||f^\infty_2(y_2)|}{{(|x-y_1|+|x-y_2|)}^{2n}}dy_1dy_2\\
  &\lesssim\dint_{2B}|f_1(y_1)|dy_1\dint_{{(2B)}^c}\frac{|f_2(y_2)|}{{|x_0-y_2|}^{2n}}dy_2 \\
   &\lesssim\dint_{2B}|f_1(y_1)|dy_1\dint_{{(2B)}^c}|f_2(y_2)|\biggl[\dint^\infty_{|x_0-y_2|}\frac{dt}{t^{2n+1}}\biggr]dy_2  \\
        &\lesssim\|f_1\|_{L^{p_1}(2B)}{|2B|}^{1-1/p_1}\dint^{\infty}_{2r}\|f_2\|_{L^{p_2}(B(x_0,t))}{|B(x_0,t)|}^{1-1/p_2}\frac{dt}{t^{2n+1}}  \\
     &\lesssim\dint^\infty_{2r}\|f_1\|_{L^{p_1}(B(x_0,t))}\|f_2\|_{L^{p_2}(B(x_0,t))}\frac{dt}{t^{n/p+1}},
\end{array}\eqno{(3.3)}$$
where $1/p=1/{p_1}+1/{p_2}.$

Thus,
$$\begin{array}{cl}II=\|T_2(f^0_1,f^\infty_2)\|_{L^p(B)} \lesssim r^{n/p}\dint^\infty_{2r}\|f_1\|_{L^{p_1}(B(x_0,t))}\|f_2\|_{L^{p_2}(B(x_0,t))}\frac{dt}{t^{n/p+1}}.\end{array}\eqno{(3.4)}$$
Similarly, we have
$$III=\|T_2(f^\infty_1,f^0_2)\|_{L^p(B)} \lesssim r^{n/p}\int^\infty_{2r}\|f_1\|_{L^{p_1}(B(x_0,t))}\|f_2\|_{L^{p_2}(B(x_0,t))}\frac{dt}{t^{n/p+1}}.$$

Moreover, similar to the estimate of (3.3), we have
 $$\begin{array}{cl}
 |T_2(f^\infty_1,f^\infty_2)(x)|
 &\lesssim\dint_{(2B)^c}\int_{(2B)^c}\frac{|f_1(y_1)||f_2(y_2)|}{{(|x_0-y_1|+|x_0-y_2|)}^{2n}}dy_1dy_2  \\
    &\lesssim\dint_{(2B)^c}\int_{(2B)^c}|f_1(y_1)||f_2(y_2)|dy_1dy_2\int^\infty_{|x_0-y_1|+|x_0-y_2|}\frac{dt}{t^{2n+1}} \\
    &\lesssim\dint^\infty_{2r}\biggl[\int_{B(x_0,t)}|f_1(y_1)|dy_1\int_{B(x_0,t)}|f_2(y_2)|dy_2\biggr]\frac{dt}{t^{2n+1}}  \\
   &\lesssim\dint^\infty_{2r}\|f\|_{L^{p_1}(B(x_0,t))}\|f\|_{L^{p_1}(B(x_0,t))}{|B(x_0,t)|}^{2-(1/p_1+1/p_2)}\frac{dt}{t^{2n+1}}    \\
    &\lesssim\dint^\infty_{2r}\|f\|_{L^{p_1}(B(x_0,t))}\|f\|_{L^{p_1}(B(x_0,t))}\frac{dt}{t^{n/p+1}}.
 \end{array}$$

 Thus,
 $$\begin{array}{cl}IV=\|T_2(f^\infty_1,f^\infty_2)\|_{L^p(B)} \lesssim r^{n/p}\dint^\infty_{2r}\|f_1\|_{L^{p_1}(B(x_0,t))}\|f_2\|_{L^{p_2}(B(x_0,t))}\frac{dt}{t^{n/p+1}}. \end{array}\eqno{(3.5)}$$

 Combining the above estimates, we obtain
  $$\begin{array}{cl}\|T_2(f_1,f_2)\|_{L^p(B)} \lesssim r^{n/p}\dint^\infty_{2r}\|f_1\|_{L^{p_1}(B(x_0,t))}\|f_2\|_{L^{p_2}(B(x_0,t))}\frac{dt}{t^{n/p+1}}.\end{array}$$

 {\bf Theorem 3.2}\;\;Let $x_0\in{\mathbb{R}^n},$ $1<p, p_1, p_2, \dots, p_m<\infty$ such that $1/p=1/p_1+1/p_2+\dots+p_m.$ If functions $\varphi,$ $\varphi_i:$ $\mathbb{R}^n\times(0,\infty)\rightarrow(0,+\infty),(i=1,2,\cdots,m)$
 satisfy the condition
 $$\begin{array}{cl}\dint_{r}^{\infty}\dfrac{\mbox{ess}\inf\limits_{t<s<\infty}\prod\limits_{i=1}^m\varphi_{i}(x_0,s)s^{n/p}}{t^{n/p+1}}dt\leq C\psi(x_0,r),\end{array}\eqno{(3.6)}$$
where constant $C>0$ doesn't depend on $r.$ Then the operator $T_m$ is bounded from the product space
$LM_{p_1,\varphi_1}^{\{x_0\}}\times LM_{p_2,\varphi_2}^{\{x_0\}}\times\dots\times LM_{p_m,\varphi_m}^{\{x_0\}}$ to $LM_{p,\psi}^{\{x_0\}}.$ Moreover,
the following inequality
 $$\|T_m(\vec{f})\|_{LM_{p,\psi}^{\{x_0\}}} \lesssim \prod\limits_{i=1}^{m}\|f_i\|_{LM_{p_i,\varphi_i}^{\{x_0\}}}.$$
holds.

{\bf Proof.} Taking $v_1(r)=\prod\limits_{i=1}^{m}\varphi_i^{-1}(x_0,r)r^{-n/p},$ $v_2(r)=\psi^{-1}(x_0,r),$
$g(r)=\prod\limits_{i=1}^{m}\|f_i\|_{L^{p_i}(B(x_0,r))}$ and $w(r)=r^{-n/p-1},$ then we have
$$ess\sup\limits_{t>0}v_2(t)\dint_{t}^{\infty}\dfrac{w(s)ds}{ess\sup\limits_{s<\tau<\infty}v_1(\tau)}<\infty.$$
 Thus, by Lemma 2.2, we have
 $$\begin{array}{cl}ess\sup\limits_{t>0}v_2(t)H_wg(t)\leq C ess\sup\limits_{t>0}v_1(t)g(t).\end{array}\eqno{(3.7)}$$

Therefore, from Theorem 3.1 and (3.7), it follows that

$$\begin{array}{cl}
&\|T_m(\vec{f})\|_{LM_{p,\psi}^{\{x_0\}}}\\
=&\sup\limits_{r>0}\psi^{-1}(x_0,r)|B(x_0,r)|^{-1/p}\|T_m(\vec{f})\|_{L^{p}(B(x_0,r))}\\
\lesssim&\sup\limits_{r>0}\psi^{-1}(x_0,r)|B(x_0,r)|^{-1/p}r^{n/p}\dint^{\infty}_{2r}\prod\limits_{i=1}^{m}\|f_i\|_{L^{p_i}(B(x_0,t))}
t^{-n/p-1}dt\\
\lesssim&\sup\limits_{r>0}\prod\limits_i^m\varphi_{i}^{-1}(x_0,r)r^{-n/p}\prod\limits_{i=1}^m\|f_i\|_{L^{p_i}(B(x_0,r))}\\
\lesssim&\sup\limits_{r>0}\prod\limits_{i=1}^m\varphi_{i}^{-1}(x_0,r)r^{-n/p_i}\|f_i\|_{L^{p_i}(B(x_0,r))}\\
=&\prod\limits_{i=1}^{m}\|f_i\|_{LM_{p_i,\varphi_i}^{\{x_0\}}}.
\end{array}$$

\section{Commutators generated by m-th Calder\'{o}n Zygmund type
singular integral operators and local Campanato functions}

{\bf Theorem 4.1}\;\; Let $x_0\in{\mathbb{R}^n},$ $1<p,$ $p_i,$ $q_i<\infty (i=1,2,\dots,m)$ such that $1/p=1/p_1+1/p_2+\dots+1/p_m+1/p_1+1/q_2+\dots+1/q_m$  and  $b_i\in LC_{q_i,\lambda_i}^{\{x_0\}}$ for $0<\lambda_i<1/n,$ $i=1,2,\cdots,m.$
Then the inequality
$$\|T_m^{\vec{b}}(\vec{f})\|_{L^p(B(x_0,r))}\lesssim \prod\limits_{i=1}^{m}\|b_i\|_{LC_{q_i,\lambda_i}^{\{x_0\}}}\;r^{n/p}
\dint^{\infty}_{2r}\biggl(1+\ln\frac{t}{r}\biggr)^{m}t^{n\sum\limits_{i=1}^{m}\lambda_i-n\sum\limits_{i=1}^{m}1/p_i-1
}\prod\limits_{i=1}^{m}\|f_i\|_{L^{p_i}(B(x_0,t))}dt$$
holds for any ball $B(x_0,r)$ and  all $f_i\in L^{p_i}_{loc}(\mathbb{R}^n),$ $i=1,2,\dots,m.$

{\bf Proof.} Without loss of generality, it is suffice for us to show that the conclusion holds for
$m=2.$

 Let $B=B(x_0,r),$ $f_1=f^0_1+f^\infty_1$ and $f_2=f^0_2+f^\infty_2,$ where $f^0_i$ and $f^\infty_i$ are as in the proof of Theorem 3.1, for $i=1.2.$ Thus, we have
$$\begin{array}{cl}&T_{2}^{(b_1,b_2)}(f_{1},f_{2})(x)\\
=&T_{2}^{(b_1,b_2)}(f_{1}^0,f_{2}^0)(x)+T_{2}^{(b_2,b_2)}(f_{1}^0,f_{2}^{\infty})(x)+T_{2}^{(b_1,b_2)}(f_{1}^{\infty},f_{2}^{0})(x)
+T_{2}^{(b_1,b_2)}(f_{1}^{\infty},f_{2}^{\infty})(x).
\end{array}$$

So,$$\begin{array}{cl}&\|T_{2}^{(b_1,b_2)}(f_{1},f_{2})\|_{L^{p}(B)}\\
&\leq\|T_{2}^{(b_1,b_2)}(f_{1}^0,f_{2}^0)\|_{L^{p}(B)}+\|T_{2}^{(b_1,b_2)}(f_{1}^0,f_{2}^{\infty})\|_{L^{p}(B)}\\
&\quad+\|T_{2}^{(b_1,b_2)}(f_{1}^{\infty},f_{2}^{0})\|_{L^{p}(B)}+\|T_{2}^{(b_1,b_2)}(f_{1}^{\infty},f_{2}^{\infty})\|_{L^{p}(B)}\\
&=:I+II+III+IV.\end{array}$$

Let us estimate $I, II, III$ and $IV$, respectively.

Since, $$\begin{array}{cl}&(b_1(x)-b_1(y))(b_2(x)-b_2(y))\\

&=(b_1(x)-(b_1)_{B})(b_2(x)-(b_2)_{B})-(b_1(x)-(b_1)_{B})(b_2(y)-(b_2)_{B})\\
&\quad-(b_1(y)-(b_1)_{B})(b_2(x)-(b_2)_{B})+(b_1(y)-(b_1)_{B})(b_2(y)-(b_2)_{B}).\end{array}\eqno{(4.1)}$$

Then,
$$\begin{array}{cl}
&\|T_{2}^{(b_1,b_2)}(f_{1}^0,f_{2}^0)\|_{L^{p}(B)}\\
=&\|(b_1-(b_1)_{B})(b_2-(b_2)_{B})T_{2}(f_{1}^0,f_{2}^0)\|_{L^{p}(B)}+\|(b_1-(b_1)_{B})T_{2}(f_{1}^0,(b_2-(b_2)_{B})f_{2}^0)\|_{L^{p}(B)}\\
&+\|(b_2-(b_2)_{B})T_{2}((b_1-(b_1)_{B})f_{1}^0,f_{2}^0)\|_{L^{p}(B)}+\|T_{2}((b_1-(b_1)_{B})f_{1}^0,(b_2-(b_2)_{B})f_{2}^0)\|_{L^{p}(B)}\\
=:&I_1+I_2+I_3+I_4.
\end{array}\eqno{(4.2)}$$

Let $1<\bar{p},\bar{q}<\infty,$ such that $1/\bar{p}=1/p_1+1/p_2$ and $1/\bar{q}=1/q_1+1/q_2.$ Then, using the H\"{o}lder's inequality and Lemma 2.3, we have
$$\begin{array}{cl}
I_1
&\lesssim\|(b_1-(b_1)_{B})(b_2-(b_2)_{B}\|_{L^{\bar{q}}(B)}\|T_{2}(f_{1}^0,f_{2}^0)\|_{L^{\bar{p}}(B)}\\
&\lesssim\|b_1-(b_1)_{B}\|_{L^{q_2}(B)}\|b_2-(b_2)_{B}\|_{L^{q_2}(B)}\|f_{1}\|_{L^{p_1}(2B)}\|f_{2}\|_{L^{p_1}(2B)}\\
&\lesssim\|b_1-(b_1)_{B}\|_{L^{q_1}(B)} \|b_2-(b_2)_{B}\|_{L^{q_2}(B)} r^{(1/p_1+1/p_2)n}\\
&\quad\times\dint^\infty_{2r}\|f_1\|_{L^{p_1}(B(x_0,t))}\|f_2\|_{L^{p_2}(B(x_0,t))}\frac{dt}{t^{(1/p_1+1/p_2)n+1}}\\

&\lesssim\|b_1\|_{LC_{q_1,\lambda_1}^{\{x_0\}}}\|b_2\|_{LC_{q_2,\lambda_2}^{\{x_0\}}}r^{n/p}\\
&\quad\times\dint^\infty_{2r}\biggl(1+\ln\dfrac{t}{r}\biggr)^2{t^{(\lambda_1+\lambda_2)n-(1/p_1+1/p_2)n-1}}\|f_1\|_{L^{p_1}(B(x_0,t))}\|f_2\|_{L^{p_2}(B(x_0,t))}dt.
\end{array}\eqno{(4.3)}$$

Let $1<\tau<\infty,$ such that $1/p=1/q_1+1/\tau.$  Then similarly to the estimate of (4.3), we have
$$\begin{array}{cl}
I_2
&\lesssim\|b_1-(b_1)_{B}\|_{L^{q_1}(B)}\|T_{2}(f_{1}^0,(b_2-(b_2)_{B})f_{2}^0)\|_{L^{\tau}(B)}\\
&\lesssim\|b_1-(b_1)_{B}\|_{L^{q_1}(B)}\|f_{1}^0\|_{L^{p_1}(\mathbb{R}^n)}\|(b_2-(b_2)_{2B})f_{2}^0)\|_{L^{s}(\mathbb{R}^n)}\\
&\lesssim\|b_1-(b_1)_{B}\|_{L^{q_1}(B)}\|b_2-(b_2)_{B}\|_{L^{q_2}(2B)}\|f_{1}\|_{L^{p_1}(2B)}\|f_{2}
\|_{L^{p_2}(2B)},\\
\end{array}\eqno{(4.4)}$$
where $1<s<\infty,$ such that $1/s=1/p_2+1/q_2=1/\tau-1/{p_1}.$

From Lemma 2.1, it is easy to see that

$$\|b_i-(b_i)_{B}\|_{L^{q_i}(B)}\leq Cr^{n/q_{i}+n\lambda_i}\|b_i\|_{LC_{q_i,\lambda_i}^{\{x_0\}}},$$
and
$$\begin{array}{cl}\|b_i-(b_i)_{B}\|_{L^{q_i}(2B)}\leq \|b_i-(b_i)_{2B}\|_{L^{q_i}(2B)}+\|(b_i)_{B}-(b_i)_{2B}\|_{L^{q_i}(2B)}\leq Cr^{n/q_{i}+n\lambda_i}\|b_i\|_{LC_{q_i,\lambda_i}^{\{x_0\}}},\end{array}\eqno{(4.5)}$$
for $i=1,2.$

Then,
$$\begin{array}{cl}
I_2
&\lesssim\|b_1\|_{LC_{q_1,\lambda_1}^{\{x_0\}}}\|b_2\|_{LC_{q_2,\lambda_2}^{\{x_0\}}}r^{n/p}\\
&\times\dint^\infty_{2r}\biggl(1+\ln\dfrac{t}{r}\biggr)^2{t^{(\lambda_1+\lambda_2)n-(1/p_1+1/p_2)n-1}}\|f_1\|_{L^{p_1}(B(x_0,t))}\|f_2\|_{L^{p_2}(B(x_0,t))}dt.
\end{array}$$

Similarly,
$$\begin{array}{cl}
I_3
&\lesssim\|b_1\|_{LC_{q_1,\lambda_1}^{\{x_0\}}}\|b_2\|_{LC_{q_2,\lambda_2}^{\{x_0\}}}r^{n/p}\\
&\times\dint^\infty_{2r}\biggl(1+\ln\dfrac{t}{r}\biggr)^2{t^{(\lambda_1+\lambda_2)n-(1/p_1+1/p_2)n-1}}\|f_1\|_{L^{p_1}(B(x_0,t))}\|f_2\|_{L^{p_2}(B(x_0,t))}dt.
\end{array}$$

Moreover, let $1<\tau_1,\tau_2<\infty,$ such that $1/\tau_{1}=1/p_1+1/q_1$ and $1/\tau_{2}=1/p_2+1/q_2.$ It is easy to see that
$1/p=1/\tau_{1}+1/\tau_{2}.$
 Then by Lemma 2.3, H\"{o}lder's inequality and (4.5), we obtain
$$\begin{array}{cl}
I_4
&\lesssim\|(b_1-(b_1)_{B})f_{1}^0\|_{L^{\tau_1}(\mathbb{R}^n)}\|(b_2-(b_2)_{B})f_{2}^0\|_{L^{\tau_2}(\mathbb{R}^n)}\\
&\lesssim\|b_1-(b_1)_{B}\|_{L^{q_1}(2B)}\|b_2-(b_2)_{B}\|_{L^{q_2}(2B)}\|f_{1}\|_{L^{p_1}(2B)}\|f_{2}
\|_{L^{p_2}(2B)}\\
&\lesssim\|b_1\|_{LC_{q_1,\lambda_1}^{\{x_0\}}}\|b_2\|_{LC_{q_2,\lambda_2}^{\{x_0\}}}r^{n/p}\\
&\times\dint^\infty_{2r}\biggl(1+\ln\dfrac{t}{r}\biggr)^2{t^{(\lambda_1+\lambda_2)n-(1/p_1+1/p_2)n-1}}\|f_1\|_{L^{p_1}(B(x_0,t))}\|f_2\|_{L^{p_2}(B(x_0,t))}dt.
\end{array}\eqno{(4.6)}$$

Therefore, combining the estimates of $I_1, I_2, I_3$ and $I_4,$ we have
$$\begin{array}{cl}I
&\lesssim\|b_1\|_{LC_{q_1,\lambda_1}^{\{x_0\}}}\|b_2\|_{LC_{q_2,\lambda_2}^{\{x_0\}}}r^{n/p}\\
&\times\dint^\infty_{2r}\biggl(1+\ln\dfrac{t}{r}\biggr)^2{t^{(\lambda_1+\lambda_2)n-(1/p_1+1/p_2)n-1}}\|f_1\|_{L^{p_1}(B(x_0,t))}\|f_2\|_{L^{p_2}(B(x_0,t))}dt.

\end{array}$$

Let us estimate $II.$

It's analogues to (4.2), we have
$$\begin{array}{cl}
&\|T_{2}^{(b_1,b_2)}(f_{1}^0,f_{2}^\infty)\|_{L^{p}(B)}\\
=&\|(b_1-(b_1)_{B})(b_2-(b_2)_{B})T_{2}(f_{1}^0,f_{2}^\infty)\|_{L^{p}(B)}+\|(b_1-(b_1)_{B})T_{2}(f_{1}^0,(b_2-(b_2)_{B})f_{2}^\infty)\|_{L^{p}(B)}\\
&+\|(b_2-(b_2)_{B})T_{2}((b_1-(b_1)_{B})f_{1}^0,f_{2}^\infty)\|_{L^{p}(B)}+\|T_{2}((b_1-(b_1)_{B})f_{1}^0,(b_2-(b_2)_{B})f_{2}^\infty)\|_{L^{p}(B)}\\
=:&II_1+II_2+I_3+II_4.
\end{array}\eqno{(4.7)}$$

Let $1<\bar{p},\bar{q}<\infty,$ such that $1/\bar{p}=1/p_1+1/p_2$ and $1/\bar{q}=1/q_1+1/q_2.$ Then, using the H\"{o}lder's inequality and (3.4), we have
$$\begin{array}{cl}
II_1
&\lesssim\|(b_1-(b_1)_{B})(b_2-(b_2)_{2B}\|_{L^{\bar{q}}(B)}\|T_{2}(f_{1}^0,f_{2}^\infty)\|_{L^{\bar{p}}(B)}\\
&\lesssim\|b_1\|_{LC_{q_1,\lambda_1}^{\{x_0\}}}\|b_2\|_{LC_{q_2,\lambda_2}^{\{x_0\}}}r^{(1/q_1+1/q_2)n+(\lambda_1+\lambda_2)n}r^{(1/p_1+1/p_2)n}\\
&\quad\times\dint^\infty_{2r}\biggl(1+\ln\dfrac{t}{r}\biggr)^2{t^{-(1/p_1+1/p_2)n-1}}\|f_1\|_{L^{p_1}(B(x_0,t))}\|f_2\|_{L^{p_2}(B(x_0,t))}dt\\
&\lesssim\|b_1\|_{LC_{q_1,\lambda_1}^{\{x_0\}}}\|b_2\|_{LC_{q_2,\lambda_2}^{\{x_0\}}}r^{n/p}\\
&\quad\times\dint^\infty_{2r}\biggl(1+\ln\dfrac{t}{r}\biggr)^2{t^{(\lambda_1+\lambda_2)n-(1/p_1+1/p_2)n-1}}\|f_1\|_{L^{p_1}(B(x_0,t))}\|f_2\|_{L^{p_2}(B(x_0,t))}dt.
\end{array}\eqno{(4.8)}$$

Moreover, using (1.2) and (3.2), we have
$$\begin{array}{cl}
&|T_{2}(f_{1}^0,(b_2-(b_2)_{B})f_{2}^\infty)(x)|\\

  &\lesssim\dint_{2B}|f_1(y_1)|dy_1\dint_{{(2B)}^c}\frac{|b_2(y_2)-(b_2)_{B}||f_2(y_2)|}{{|x_0-y_2|}^{2n}}dy_2. \\
\end{array}$$

It's obvious that
$$\begin{array}{cl}\dint_{2B}|f_1(y_1)|dy_1\lesssim \|f_1\|_{L^{p_1}(2B)}{|2B|}^{1-1/p_1},\end{array}\eqno{(4.9)}$$
and

$$\begin{array}{cl}
&\dint_{{(2B)}^c}\frac{|b_2(y_2)-(b_2)_{B}||f_2(y_2)|}{{|x_0-y_2|}^{2n}}dy_2 \\
   &\lesssim\dint_{{(2B)}^c}|b_2(y_2)-(b_2)_{B}||f_2(y_2)|\biggl[\dint^\infty_{|x_0-y_2|}\frac{dt}{t^{2n+1}}\biggr]dy_2  \\
         &\lesssim \dint^{\infty}_{2r}
   \|b_2(y_2)-(b_2)_{B(x_0,t)}\|_{L^{q_2}(B(x_0,t))}\|f_2\|_{L^{p_2}(B(x_0,t))}|B(x_0,t)|^{1-(1/p_2+1/q_2)}\frac{dt}{t^{2n+1}}\\
    &\quad+\dint^{\infty}_{2r}|(b_2)_{B(x_0,t)}-(b_2)_{B(x_0,r)}
    |\|f_2\|_{L^{p_2}(B(x_0,t))}|B(x_0,t)|^{1-1/p_2}\frac{dt}{t^{2n+1}}\\
      &\lesssim \|b_2\|_{LC_{q_2,\lambda_2}^{\{x_0\}}}
   \dint^{\infty}_{2r}|B(x_0,t)|^{1/q_2+\lambda_2}\|f_2\|_{L^{p_2}(B(x_0,t))}|B(x_0,t)|^{1-(1/p_2+1/q_2)}\frac{dt}{t^{2n+1}}\\
    &\quad+\|b_2\|_{LC_{q_2,\lambda_2}^{\{x_0\}}}
\dint^{\infty}_{2r}\biggl(1+\ln\dfrac{t}{r}\biggr)|B(x_0,t)|^{\lambda_2}
    \|f_2\|_{L^{p_2}(B(x_0,t))}|B(x_0,t)|^{1-1/p_2}\frac{dt}{t^{2n+1}}\\
   &\lesssim\|b_2\|_{LC_{q_2,\lambda_2}^{\{x_0\}}}
    \dint^{\infty}_{2r}\biggl(1+\ln\dfrac{t}{r}\biggr)^2t^{-n+n\lambda_2-n/p_2-1}\|f_2\|_{L^{p_2}(B(x_0,t))}dt.\\
\end{array}\eqno{(4.10)}$$

Therefore, from (4.9) and (4.10), it follows that
$$\begin{array}{cl}
&|T_{2}(f_{1}^0,(b_2-(b_2)_{B})f_{2}^\infty)(x)|\\
&\lesssim\|b_2\|_{LC_{q_2,\lambda_2}^{\{x_0\}}}\|f_1\|_{L^{p_1}(2B)}{|2B|}^{1-1/p_1}
    \dint^{\infty}_{2r}\biggl(1+\ln\dfrac{t}{r}\biggr)^2t^{-n+n\lambda_2-n/p_2-1}\|f_2\|_{L^{p_2}(B(x_0,t))}dt\\
&\lesssim\|b_2\|_{LC_{q_2,\lambda_2}^{\{x_0\}}} \dint^{\infty}_{2r}\biggl(1+\ln\dfrac{t}{r}\biggr)^2t^{n\lambda_2-(1/p_1+1/p_2)n-1}\|f_1\|_{L^{p_1}(B(x_0,t))}\|f_2\|_{L^{p_2}(B(x_0,t))}dt.\\
\end{array}$$

Thus, let $1<\tau<\infty,$ such that $1/p=1/q_1+1/\tau$, then similarly to the estimate of (4.3), we have
$$\begin{array}{cl}
II_2&=\|(b_1-(b_1)_{B})T_{2}(f_{1}^0,(b_2-(b_2)_{B})f_{2}^\infty)\|_{L^p(B)}\\
&\lesssim\|b_1-(b_1)_{B}\|_{L^{q_1}(B)}\|T_{2}(f_{1}^0,(b_2-(b_2)_{B})f_{2}^\infty)\|_{L^{\tau}(B)}\\
&\lesssim\|b_1\|_{LC_{q_1,\lambda_1}^{\{x_0\}}}\|b_2\|_{LC_{q_2,\lambda_2}^{\{x_0\}}}|B|^{\lambda_1+1/{q_1}+1/\tau}\\
&\quad\times\dint^{\infty}_{2r}\biggl(1+\ln\dfrac{t}{r}\biggr)^2t^{n\lambda_2-(1/p_1+1/p_2)n-1}\|f_1\|_{L^{p_1}(B(x_0,t))}\|f_2\|_{L^{p_2}(B(x_0,t))}dt\\
&\lesssim\|b_1\|_{LC_{q_1,\lambda_1}^{\{x_0\}}}\|b_2\|_{LC_{q_2,\lambda_2}^{\{x_0\}}}r^{n/p}\\
&\quad\times\dint^\infty_{2r}\biggl(1+\ln\dfrac{t}{r}\biggr)^2{t^{(\lambda_1+\lambda_2)n-(1/p_1+1/p_2)n-1}}\|f_1\|_{L^{p_1}(B(x_0,t))}\|f_2\|_{L^{p_2}(B(x_0,t))}dt.
\end{array}\eqno{(4.11)}$$

Similarly, we have
$$\begin{array}{cl}
II_3&\lesssim\|b_1\|_{LC_{q_1,\lambda_1}^{\{x_0\}}}\|b_2\|_{LC_{q_2,\lambda_2}^{\{x_0\}}}r^{n/p}\\
&\quad\times\dint^\infty_{2r}\biggl(1+\ln\dfrac{t}{r}\biggr)^2{t^{(\lambda_1+\lambda_2)n-(1/p_1+1/p_2)n-1}}\|f_1\|_{L^{p_1}(B(x_0,t))}\|f_2\|_{L^{p_2}(B(x_0,t))}dt.
\end{array}$$

Let us estimate $II_4.$

Since,
$$\begin{array}{cl}
&|T_{2}((b_1-(b_1)_{B})f_{1}^0,(b_2-(b_2)_{B})f_{2}^\infty)(x)|\\
  &\lesssim\dint_{2B}|b_1(y_1)-(b_1)_{B}||f_1(y_1)|dy_1\dint_{{(2B)}^c}\frac{|b_2(y_2)-(b_2)_{B}||f_2(y_2)|}{{|x_0-y_2|}^{2n}}dy_2, \\
\end{array}$$
and
$$\begin{array}{cl}&\dint_{2B}|b_1(y_1)-(b_1)_{B}||f_1(y_1)|dy_1
  \lesssim\|b_1\|_{LC_{q_1,\lambda_1}^{\{x_0\}}}|B|^{\lambda_1+1-1/p_1}\|f_1\|_{L^{p_1}(2B)}.
\end{array}\eqno{(4.12)}$$

Then, by (4.10) and (4.12), we have
$$\begin{array}{cl}
&|T_{2}((b_1-(b_1)_{B})f_{1}^0,(b_2-(b_2)_{B})f_{2}^\infty)(x)|\\

&\lesssim\|b_1\|_{LC_{q_1,\lambda_1}^{\{x_0\}}}\|b_2\|_{LC_{q_2,\lambda_2}^{\{x_0\}}}
    \dint^{\infty}_{2r}\biggl(1+\ln\dfrac{t}{r}\biggr)^2t^{n(\lambda_1+\lambda_2)-n(1/p_1+/p_2)-1}\|f_1\|_{L^{p_1}(B(x_0,t))}\|f_2\|_{L^{p_2}(B(x_0,t))}dt.
\end{array}$$

Therefore,
$$\begin{array}{cl}
II_4&=\|T_{2}((b_1-(b_1)_{B})f_{1}^0,(b_2-(b_2)_{B})f_{2}^\infty)\|_{L^p(B)}\\
&\lesssim\|b_1\|_{LC_{q_1,\lambda_1}^{\{x_0\}}}\|b_2\|_{LC_{q_2,\lambda_2}^{\{x_0\}}}r^{n/p}\\
&\quad\times\dint^\infty_{2r}\biggl(1+\ln\dfrac{t}{r}\biggr)^2{t^{(\lambda_1+\lambda_2)n-(1/p_1+1/p_2)n-1}}\|f_1\|_{L^{p_1}(B(x_0,t))}\|f_2\|_{L^{p_2}(B(x_0,t))}dt.
\end{array}$$

Combining the estimates of $II_1-II_4,$ we have
$$\begin{array}{cl}
II
&\lesssim\|b_1\|_{LC_{q_1,\lambda_1}^{\{x_0\}}}\|b_2\|_{LC_{q_2,\lambda_2}^{\{x_0\}}}r^{n/p}\\
&\quad\times\dint^\infty_{2r}\biggl(1+\ln\dfrac{t}{r}\biggr)^2{t^{(\lambda_1+\lambda_2)n-(1/p_1+1/p_2)n-1}}\|f_1\|_{L^{p_1}(B(x_0,t))}\|f_2\|_{L^{p_2}(B(x_0,t))}dt.
\end{array}$$

Similarly,
$$\begin{array}{cl}
III&\lesssim\|b_1\|_{LC_{q_1,\lambda_1}^{\{x_0\}}}\|b_2\|_{LC_{q_2,\lambda_2}^{\{x_0\}}}r^{n/p}\\
&\quad\times\dint^\infty_{2r}\biggl(1+\ln\dfrac{t}{r}\biggr)^2{t^{(\lambda_1+\lambda_2)n-(1/p_1+1/p_2)n-1}}\|f_1\|_{L^{p_1}(B(x_0,t))}\|f_2\|_{L^{p_2}(B(x_0,t))}dt.
\end{array}$$

For $IV,$ we have
$$\begin{array}{cl}
&\|T_{2}^{(b_1,b_2)}(f_{1}^\infty,f_{2}^{\infty})\|_{L^{p}(B)}\\
\leq&\|(b_1-(b_1)_{B})(b_2-(b_2)_{B})T_{2}(f_{1}^\infty,f_{2}^{\infty})\|_{L^{p}(B)}
+\|(b_1-(b_1)_{B})T_{2}(f_{1}^\infty,(b_2-(b_2)_{B})f_{2}^{\infty})\|_{L^{p}(B)}\\
&\quad+\|(b_2-(b_2)_{B})T_{2}((b_1-(b_1)_{B})f_{1}^\infty,f_{2}^{\infty})\|_{L^{p}(B)}+\|T_{2}((b_1-(b_1)_{B})f_{1}^\infty,(b_2-(b_2)_{B})f_{2}^{\infty})\|_{L^{p}(B)}\\
=:&IV_1+IV_2+IV_3+IV_4.\end{array}$$

Let us estimate $IV_1,$ $IV_2,$ $IV_3$ and $IV_4,$ respectively.

Let $1<\tau<\infty,$ such that $1/p=1/q_1+1/q_2+1/\tau.$ Then, from H\"{o}lder's inequality and (3.5), we get
$$\begin{array}{cl}
IV_1&\lesssim\|b_1-(b_1)_{B}\|_{L^{q_1}(B)}\|b_2-(b_2)_{B}\|_{L^{q_2}(B)}\|T_{2}(f_{1}^\infty,f_{2}^{\infty})\|_{L^{\tau}(B)}\\
&\lesssim\|b_1\|_{LC_{q_1,\lambda_1}^{\{x_0\}}}\|b_2\|_{LC_{q_2,\lambda_2}^{\{x_0\}}}|B|^{(\lambda_1+\lambda_2)+(1/q_1+1/q_2)+1/\tau}\\
&\quad\times\dint^\infty_{2r}\|f_1\|_{L^{p_1}(B(x_0,t))}\|f_2\|_{L^{p_2}(B(x_0,t))}t^{-n(/p_1+1/p_2)-1}dt\\
&\lesssim\|b_1\|_{LC_{q_1,\lambda_1}^{\{x_0\}}}\|b_2\|_{LC_{q_2,\lambda_2}^{\{x_0\}}}r^{n/p}\\
&\quad\times\dint^\infty_{2r}\biggl(1+\ln\dfrac{t}{r}\biggr)^2{t^{(\lambda_1+\lambda_2)n-(1/p_1+1/p_2)n-1}}\|f_1\|_{L^{p_1}(B(x_0,t))}\|f_2\|_{L^{p_2}(B(x_0,t))}dt.
\end{array}$$

Moreover, by (1.2) and (3.2), we have
$$\begin{array}{cl}
&|T_{2}(f_{1}^\infty,(b_2-(b_2)_{B})f_{2}^\infty)(x)|\\
&\lesssim\dint_{{(2B)}^c}\dint_{{(2B)}^c}\dfrac{|b_2(y_2)-(b_2)_{B}||f_1(y_1)||f_2(y_2)|}{{(|x_0-y_1|+|x_0-y_2|)}^{2n}}dy_1dy_2\\
 &\lesssim\dint_{{(2B)}^c}\dint_{{(2B)}^c}|f_1(y_1)||b_2(y_2)-(b_2)_{B}||f_2(y_2)|\biggl[\dint^\infty_{|x_0-y_1|+|x_0-y_2|}\frac{dt}{t^{2n+1}}\biggr]dy_1dy_2 \\
    &\lesssim \dint^{\infty}_{2r}\biggr[\dint_{B(x_0,t)}|f_1(y_1)|dy_1\biggr]\biggr[\dint_{B(x_0,t)}|b_2(y_2)-(b_2)_{B}|
    |f_2(y_2)|dy_2\biggr]\frac{dt}{t^{2n+1}}.\\
 \end{array}$$

 Since,

$$\begin{array}{cl}\dint_{B(x_0,t)}|f_1(y_1)|dy_1\lesssim\|f_1\|_{L^{p_1}(B(x_0,t))}t^{n(1-1/p_1)}, \end{array}$$
and
$$\begin{array}{cl}
   &\dint_{B(x_0,t)}|b_2(y_2)-(b_2)_{B}||f_2(y_2)|  \\
         &\lesssim\|b_2-(b_2)_{B(x_0,t)}\|_{L^{q_2}(B(x_0,t))}\|f_2\|_{L^{p_2}}|B(x_0,t)|^{1-(1/p_2+1/q_2)}\\
    &\quad+|(b_2)_{B(x_0,t)}-(b_2)_{B(x_0,r)}
    |\|f_2\|_{L^{p_2}}|B(x_0,t)|^{1-1/p_2}\\
           &\lesssim \|b_2\|_{LC_{q_2,\lambda_2}^{\{x_0\}}}|B(x_0,t)|^{1/q_2+\lambda_2}\|f_2\|_{L^{p_2}}|B(x_0,t)|^{1-(1/p_2+1/q_2)}\\
      &\quad+\|b_2\|_{LC_{q_2,\lambda_2}^{\{x_0\}}}\biggl(1+\ln\dfrac{t}{r}\biggr)|B(x_0,t)|^{\lambda_2} \|f_2\|_{L^{p_2}}|B(x_0,t)|^{1-1/p_2}\\
   &\lesssim\|b_2\|_{LC_{q_2,\lambda_2}^{\{x_0\}}}
   \biggl(1+\ln\dfrac{t}{r}\biggr)^2t^{n\lambda_2-n/p_2+n}\|f_2\|_{L^{p_2}(B(x_0,t))}.\\
\end{array}$$

Then,
$$\begin{array}{cl}
&|T_{2}(f_{1}^\infty,(b_2-(b_2)_{B})f_{2}^\infty)(x)|\\
  &\lesssim\|b_2\|_{LC_{q_2,\lambda_2}^{\{x_0\}}}\dint^{\infty}_{2r}\biggl(1+\ln\dfrac{t}{r}\biggr)^2t^{n\lambda_2-(1/p_1+1/p_2)n-1}\|f_1\|_{L^{p_1}(B(x_0,t))}\|f_2\|_{L^{p_2}(B(x_0,t))}dt.\\
\end{array}\eqno{(4.13)}$$

Let $1<\tau<\infty,$ such that $1/p=1/q_1+1/\tau.$ Then, from H\"{o}lder's inequality and (4.13), we have
$$\begin{array}{cl}
IV_2
&\lesssim\|b_1-(b_1)_{B}\|_{L^{q_1}(B)}\|T_{2}(f_{1}^\infty,(b_2-(b_2)_{B})f_{2}^\infty)\|_{L^{\tau}(B)}\\
&\lesssim\|b_1\|_{LC_{q_1,\lambda_1}^{\{x_0\}}}\|b_2\|_{LC_{q_2,\lambda_2}^{\{x_0\}}}r^{n/p}\\
&\quad\times\dint^\infty_{2r}\biggl(1+\ln\dfrac{t}{r}\biggr)^2{t^{(\lambda_1+\lambda_2)n-(1/p_1+1/p_2)n-1}}\|f_1\|_{L^{p_1}(B(x_0,t))}\|f_2\|_{L^{p_2}(B(x_0,t))}dt.
\end{array}$$

Similarly,
$$\begin{array}{cl}
IV_3
&\lesssim\|b_1\|_{LC_{q_1,\lambda_1}^{\{x_0\}}}\|b_2\|_{LC_{q_2,\lambda_2}^{\{x_0\}}}r^{n/p}\\
&\quad\times\dint^\infty_{2r}\biggl(1+\ln\dfrac{t}{r}\biggr)^2{t^{(\lambda_1+\lambda_2)n-(1/p_1+1/p_2)n-1}}\|f_1\|_{L^{p_1}(B(x_0,t))}\|f_2\|_{L^{p_2}(B(x_0,t))}dt.
\end{array}$$

Similar to the estimate of (4.13), we have
$$\begin{array}{cl}
&|T_{2}((b_1-(b_1)_{B})f_{1}^\infty,(b_2-(b_2)_{B})f_{2}^\infty)(x)|\\
 &\lesssim\dint_{{(2B)}^c}\dint_{{(2B)}^c}|b_1(y_1)-(b_1)_{B}||b_2(y_2)-(b_2)_{B}||f_1(y_1)||f_2(y_2)|\biggl[\dint^\infty_{|x_0-y_1|+|x_0-y_2|}\frac{dt}{t^{2n+1}}\biggr]dy_1dy_2 \\    &\lesssim \dint^{\infty}_{2r}\biggr[\dint_{B(x_0,t)}|b_1(y_1)-(b_1)_{B}||f_1(y_1)|dy_1\biggr]\biggr[\dint_{B(x_0,t)}|b_2(y_2)-(b_2)_{B}||f_2(y_2)|dy_2\biggr]\frac{dt}{t^{2n+1}}\\
 &\lesssim \|b_1\|_{LC_{q_1,\lambda_1}^{\{x_0\}}} \|b_2\|_{LC_{q_2,\lambda_2}^{\{x_0\}}}\dint^{\infty}_{2r}
   \biggl(1+\ln\dfrac{t}{r}\biggr)^2t^{n(\lambda_1+\lambda_2)-n(1/p_1+1/p_2)-1}\|f_1\|_{L^{p_1}(B(x_0,t))}\|f_2\|_{L^{p_2}(B(x_0,t))}
dt.\\
 \end{array}$$

Thus,
$$\begin{array}{cl}
IV_4
&\lesssim\|b_1\|_{LC_{q_1,\lambda_1}^{\{x_0\}}}\|b_2\|_{LC_{q_2,\lambda_2}^{\{x_0\}}}r^{n/p}\\
&\quad\times\dint^\infty_{2r}\biggl(1+\ln\dfrac{t}{r}\biggr)^2{t^{(\lambda_1+\lambda_2)n-(1/p_1+1/p_2)n-1}}\|f_1\|_{L^{p_1}(B(x_0,t))}\|f_2\|_{L^{p_2}(B(x_0,t))}dt.
\end{array}$$

Then, from the estimates of $IV_1-IV_4,$ we deduce that
$$\begin{array}{cl}
IV&\lesssim\|b_1\|_{LC_{q_1,\lambda_1}^{\{x_0\}}}\|b_2\|_{LC_{q_2,\lambda_2}^{\{x_0\}}}r^{n/p}\\
&\quad\times\dint^\infty_{2r}\biggl(1+\ln\dfrac{t}{r}\biggr)^2{t^{(\lambda_1+\lambda_2)n-(1/p_1+1/p_2)n-1}}\|f_1\|_{L^{p_1}(B(x_0,t))}\|f_2\|_{L^{p_2}(B(x_0,t))}dt.
\end{array}$$

So, combining the estimates for $I, II, III$ and $IV$, we have
$$\begin{array}{cl}&\|T_{2}^{(b_1,b_2)}(f_{1},f_{2})\|_{L^{p}(B)}\\
&\lesssim\|b_1\|_{LC_{q_1,\lambda_1}^{\{x_0\}}}\|b_2\|_{LC_{q_2,\lambda_2}^{\{x_0\}}}r^{n/p}\\
&\quad\times\dint^\infty_{2r}\biggl(1+\ln\dfrac{t}{r}\biggr)^2{t^{(\lambda_1+\lambda_2)n-(1/p_1+1/p_2)n-1}}\|f_1\|_{L^{p_1}(B(x_0,t))}\|f_2\|_{L^{p_2}(B(x_0,t))}dt.
\end{array}$$

Therefore, we complete the proof of Theorem 4.1.

 {\bf Theorem 4.2}\;\;  Let $x_0\in{\mathbb{R}^n},$ $1<p, p_i, q_i<\infty,$ for $i=1, 2,\dots,m$ such that $1/p=1/p_1+1/p_2+\dots+1/p_n+1/p_1+1/q_2+\dots+1/q_n.$  Suppose that $0<\lambda_i<1/n$ such that $b_i\in LC_{q_i,\lambda_i}^{\{x_0\}},$ for $0<\lambda_i<1/n,$ $i=1,2,\cdots,m.$
 If functions $\varphi,$ $\varphi_i:$ $\mathbb{R}^n\times(0,\infty)\rightarrow(0,+\infty),(i=1,2,\cdots,m)$
 satisfy the condition
  $$\dint_{r}^{\infty}\biggl(1+\ln\frac{t}{r}\biggr)^{m}t^{n\sum\limits_{i=1}^{m}\lambda_i-n\sum\limits_{i=1}^{m}1/p_i-1}\mbox{ess}\inf\limits_{t<s<\infty}\prod\limits_{i=1}^m\varphi_{i}(x_0,s)s^{n/p_i}dt\leq C\psi(x_0,r),$$
where constant $C>0$ doesn't depend on $r.$ Then the operator $T_m^{\vec{b}}$ is bounded from product space
$LM_{p_1,\varphi_1}^{\{x_0\}}\times LM_{p_2,\varphi_2}^{\{x_0\}}\times\dots\times LM_{p_m,\varphi_m}^{\{x_0\}}$ to $LM_{p,\psi}^{\{x_0\}}.$ Moreover,
the inequality
 $$\|T_m^{\vec{b}}(\vec{f})\|_{LM_{p,\psi}^{\{x_0\}}}\lesssim\prod\limits_{i=1}^{m}\|b_i\|_{LC_{q_i,\lambda_i}^{\{x_0\}}}\prod\limits_{i=1}^{m}\|f_i\|_{LM_{p_i,\varphi_i}^{\{x_0\}}}.$$
holds.

{\bf Proof.} Taking $v_1(t)=\prod\limits_{i=1}^{m}\varphi_i^{-1}(x_0,t)t^{-n/p_i},$ $v_2(t)=\psi^{-1}(x_0,t),$
$g(t)=\prod\limits_{i=1}^{m}\|f_i\|_{L^{p_i}(B(x_0,t))}$ and $w(t)=(1+\ln\frac{t}{r})^{m}t^{n\sum\limits_{i=1}^{m}\lambda_i-n\sum\limits_{i=1}^{m}1/p_i-1},$
then we have
 $$ess\sup\limits_{t>0}v_2(t)\dint_{t}^{\infty}\dfrac{w(s)ds}{ess\sup\limits_{s<\tau<\infty}v_1(\tau)}<\infty.$$
 Thus, by Lemma 2.2, we have
 $$ess\sup\limits_{t>0}v_2(t)H_wg(t)\leq C ess\sup\limits_{t>0}v_1(t)g(t).$$

So,
$$\begin{array}{cl}
&\|T_m^{\vec{b}}(\vec{f})\|_{LM_{p,\psi}^{\{x_0\}}}\\
=&\sup\limits_{r>0}\psi^{-1}(x_0,r)|B(x_0,r)|^{-1/p}\|T_m(\vec{f})\|_{L^{p}(B(x_0,r))}\\
\lesssim&\prod\limits_{i=1}^{m}\|b_i\|_{LC_{q_i,\lambda_i}^{\{x_0\}}}\sup\limits_{r>0}\psi^{-1}(x_0,r)
\dint^{\infty}_{2r}\biggl(1+\ln\frac{t}{r}\biggr)^mt^{n\sum\limits_{i=1}^{m}\lambda_i-n\sum\limits_{i=1}^{m}1/p_i-1}\prod\limits_{i=1}^{m}\|f_i\|_{L^{p_i}(B(x_0,t))}dt\\
\lesssim&\prod\limits_{i=1}^{m}\|b_i\|_{LC_{q_i,\lambda_i}^{\{x_0\}}}\sup\limits_{r>0}\prod\limits_{i=1}^m\varphi_{i}^{-1}(x_0,r)r^{-n/p_i}\|f_i\|_{L^{p_i}(B(x_0,r))}\\
=&\prod\limits_{i=1}^{m}\|f_i\|_{LM_{p_i,\varphi_i}^{\{x_0\}}}. \end{array}$$
Thus we complete the proof of Theorem 4.2.

\noindent{\bf Acknowledgments}\;\;{This work is supported by the National Natural Science Foundation of China (11161042, 11471050)}

\end{document}